\documentclass[a4paper,11pt]{article}
\usepackage{amsmath,amssymb}
\usepackage{parskip}
\usepackage{hyperref}

\newtheorem{thm}{Theorem}
\newtheorem{prop}[thm]{Proposition}
\newtheorem{lemma}[thm]{Lemma}

\newtheorem{qn}{Question}

\begin{document}
\title{Torsion growth in finitely presented pro-$p$ groups} 
\author{Nikolay Nikolov}
\maketitle 

\begin{abstract} We prove that torsion in the abelianizations of open normal subgroups in finitely presented pro-$p$ groups can grow arbitrarily fast. By way of contrast in $p$-adic analytic pro-$p$ groups the torsion growth is at most polynomial. \end{abstract}

\section{Main result and motivation.} 

For a finitely generated abelian group $A$ we denote by $t(A)$ the size of the torsion subgroup of $A$. Let $p$ be a prime number and let $t_p(A)$ be the size of the Sylow $p$-subgroup of $A$, that is the largest prime power  $p^k$ which divides $t(A)$.

Let $G$ be a finitely generated residually finite group. We are interested in the growth of $t(H^{ab})$ as a function of $|G:H|$ where $H$ ranges over normal subgroups of finite index in $G$. It turns out that this growth can be arbitrary fast even when $G$ is a finitely generated solvable group \cite{KKN}. On the other hand when $G$ is finitely presented then $t(H^{ab})$ grows at most exponentially in terms of $|G:H|$, see \cite{AGN} Lemma 27. Of particular interest is the study of $t(H^{ab})$ when $G$ is an arithmetic hyperbolic 3-manifold group and $H$ ranges over its congruence subgroups, see \cite{BSV}. Proving lower bounds on $t(H^{ab})$ appears to be difficult and the following is still open.

\begin{qn} Is there a finitely presented residually finite group $G$ with a chain $G>G_1>G_2> \cdots $ of normal subgroups of finite index in $G$ such that $\cap_{i=1}^\infty G_i=\{1\}$ and  
\[ \liminf_{i \rightarrow \infty}  \frac{\log t(G_i^{ab})}{|G:G_i|}>0?\] 
\end{qn}

It may be more feasible to study $t_p(G_i^{ab})$ for single prime $p$. Motivated by this we consider torsion growth of pro-$p$ groups. 

Let $A$ be a finitely generated abelian profinite group. We extend the definition of $t(A)$ to be again the size of the torsion subgroup of $A$. However $t(A)$ can be infinite even when $A$ is procyclic, for example consider $A=\prod_{i=1}^ \infty \mathbb Z/p_i \mathbb Z$, where $p_1<p_2< \cdots $ is an infinite sequence of distinct prime numbers. Henceforth we restrict the discussion to finitely generated virtually pro-$p$ groups $\Gamma$, in which case $t(\Gamma^{ab})$ is always finite. 

The first question to answer is whether there is an asymptotic upper bound for the torsion growth of a finitely presented pro-$p$ group $\Gamma$ similar to the situation in discrete groups. The relations defining $\Gamma$ and its open subgroups are now pro-$p$ elements which may not be limits of words of bounded length in the generating set of $\Gamma$. This is a key tool in estimating torsion in the discrete case and therefore we should not expect the existence of a general upper bound for the torsion growth of $\Gamma$.
The main result of this note is to make this heuristic observation more precise.

\begin{thm} \label{main} Let $f: \mathbb N \rightarrow \mathbb N$ be any function. There exists a pro-$p$ group $\Gamma$ presented by two generators and two relations, and an infinite chain of open normal subgroups $\Gamma>\Gamma_1>\Gamma_2> \cdots $ of $\Gamma$ such that \[ \log_p t(\Gamma_i^{ab})> f(\log_p|\Gamma:\Gamma_i|) \quad \textrm{for all} \  i \in \mathbb N.\]
\end{thm}

We remark that our argument does not give that the chain $\Gamma_i$ has trivial intersection and we conjecture that the conclusion of the theorem still holds under this extra condition.  

It will be interesting to find examples of solvable finitely presented pro-$p$ groups $\Gamma$ with arbitrary fast torsion growth. The groups $\Gamma$ constructed in our proof of Theorem \ref{main} are rather big: each has a homomorphic image with a presentation of positive \emph{$p$-deficiency} as defined by Schlage-Puchta in \cite{JCP}. 

 We contrast Theorem \ref{main} with torsion growth in $p$-adic analytic groups. These groups have relatively rigid stricture tightly controlled by their Lie algebra via the Lazard correspondence. Using this and some elementary estimates we show that the growth of torsion in $p$-adic analytic groups is at most polynomial.

For a $p$-adic analytic group $G$ we denote by $\dim G$ the dimension of $G$ as an analytic group, that is the dimension of the $\mathbb Q_p$-Lie algebra of $G$.  When $G$ is an abelian pro-$p$ group then $\dim G$ is the same as the torsion free rank of $G$.

\begin{thm} \label{padic} Let $G$ be a $p$-adic analytic group. There is a positive integer $b \in \mathbb N$ such that $t(H^{ab}) \leq b |G:H|^{2\dim G}$ for every open normal subgroup $H$ of $G$. 
\end{thm}

Before we proceed to the proofs of Theorems \ref{main} and \ref{padic} we collect several preliminary results, the discrete analogues of which may be useful towards Question 1. 
\section{Preliminary results}
For a group $G$ and an integer $n$ we denote by $G[n]:=G/G^n $ the largest homomorphic image of $G$ with exponent dividing $n$. By $G^{(n)}$ we will denote the direct product of $n$ copies of $G$ and $\mathbb Z_p$ is the (additive) group of the $p$-adic integers. We begin with the following basic result. 
\begin{lemma} \label{L1} Let $A$ and $B$ be two abelian groups.

(1) Suppose $A$ and $B$ are finitely generated abelian pro-$p$ groups and let $t(A)=p^k$. Let $n>k$ be an integer and suppose $A[p^n] \simeq B[p^n]$. Then $t(B) \geq t(A)$.

(2) Suppose $A$ and $B$ are finitely generated abstract groups and let $n=t(A)m$, where $m>1$ is an integer. Suppose that $A[n] \simeq B[n]$, then $t(B) \geq t(A)$. 
\end{lemma}
\textbf{Proof:} We give the proof of (1) and the proof of (2) (which we don't need in this paper) is very similar.

Let $A=T \oplus (\mathbb Z_p)^{(a)}$ and $B=S \oplus (\mathbb Z_p)^{(b)}$ where $T$ and $S$ are finite abelian $p$-groups and $a, b$ are nonnegative integers.  Since $|T|=t(A)$ divides $p^n$ we have 
\[ A[p^n]=T \oplus (\mathbb Z/p^n \mathbb Z)^{(a)} \simeq S[p^n] \oplus (\mathbb Z/p^n \mathbb Z)^{(b)}=B[p^n].\]

We claim that $a \geq b$. Suppose to the contrary that $a<b$. We have $p^k A[p^n]=(p^k\mathbb Z/ p^n \mathbb Z)^{(a)}$ while $p^kB[p^n]$ has size at least $(p^k\mathbb Z/p^ n\mathbb Z)^{(b)}$ and so $|p^kB[p^n]|>|p^kA[p^n]|$, a contradiction proving the claim.  Therefore 
\[ |S[p^n]| p^{nb}=|B[p^n]|= |A[p^n]|=p^k p^{na} \geq p^{k} p^{nb} \] giving that $|S[p^n]| \geq p^k=|T|$ and in particular $t(B)=|S| \geq |T|=t(A)$. The Lemma follows.
$\square$

We will frequently use the following elementary result.
\begin{prop}\label{el} Let $A$ be an abelian group with a subgroup $B \leq A$ of finite index. Then $t(A) \leq t(B)|A:B|$.
\end{prop}
\textbf{Proof:} Let $T$ be the  torsion subgroup of $A$. Then $T \cap B$ is the torsion subgroup of $B$ and so $|T \cap B| = t(B)$. Together with $|T: (B \cap T)| \leq |A:B|$ this gives $|T| \leq t(B) |A:B|$ and the proposition follows. $\square$
\begin{lemma} \label{module} Let $G$ be a finite group and let $p$ be a prime. Let $M$ be a finitely generated $\mathbb Z G$-module and let $m_1, \ldots, m_d \in M$ generate a submodule of infinite index in $M$. Assume that $M$ is torsion free as a $\mathbb Z$-module and further that $M \otimes_{\mathbb Z} \mathbb Q$ is $d$-generated as $\mathbb Q G$-module. There exist elements $h_1, \ldots, h_d \in M$ such that if we define $K_n$ to be the $\mathbb Z G$-submodule of $M$ generated by $m_1+p^nh_1, \ldots, m_d+p^nh_d$ then
\[ t_p( M/K_n) \rightarrow \infty \quad \mathrm{as} \ n \rightarrow \infty.\]
\end{lemma}

\textbf{Proof:} Let $U= \sum_{i=1}^d \mathbb Z  G m_i$ be the submodule of $M$ generated by the $m_1, \ldots, m_d$ and put $\mathcal U= \mathbb Q \otimes_{\mathbb Z} U$ and $\mathcal M=\mathbb Q  \otimes_{\mathbb Z} M$. We consider $U$ and $M$ as $\mathbb Z G$-submodules of $\mathcal U$ and $\mathcal M$ respectively via the injection $m \mapsto 1 \otimes m$ for all $m \in M$.
Since $\mathbb Q G$-modules are completely reducible we can write $\mathcal M= \mathcal U \oplus \mathcal V$ for some $\mathbb Q G$-module $\mathcal V$. Morever $\mathcal V \not =\{0\}$ since $U$ has infinite index in $M$ and hence $\dim_ \mathbb Q \mathcal U < \dim_{\mathbb Q} \mathcal M$. Let $V= M \cap \mathcal V$, then $V$ is a $\mathbb Z $-lattice in the $\mathbb Q $-vector space $\mathcal V$. Note that $U \cap V \leq \mathcal U \cap \mathcal V=\{0\}$ and hence $M$ contains the direct sum $U \oplus V$ as a $\mathbb Z G$-submodule of finite index.

Consider the  $\mathbb Q G$-module homomorphism  $f: (\mathbb Q G)^{(d)} \rightarrow \mathcal M$ given by \\ $f(r_1, \ldots, r_d)= \sum_{i=1}^d r_im_i$ for all $r_i \in \mathbb Q G$. By definition $Im(f)=\mathcal U$. Let $\mathcal L= \ker f$ and $L=\mathcal L \cap (\mathbb Z G)^{(d)}$.  

By the complete reducibility of $\mathbb Q G$ again we have that $(\mathbb Q G)^{(d)} = \mathcal L \oplus \mathcal X$ with $\mathcal X \simeq \mathcal U$ as $\mathbb Q G$-modules. Since $\mathcal M$ is $d$-denerated as $\mathbb Q G$-module it follows that $(\mathbb Q G)^{(d)} \simeq \mathcal U \oplus \mathcal L$  has a direct summand isomorphic to $\mathcal M=\mathcal U \oplus \mathcal V$. By comparing multiplicites of the  simple composition factors in $\mathcal V$ and $\mathcal L$ we deduce that $\mathcal L$ has a direct summand isomorphic  to $\mathcal V$. In particular there is a surjective $\mathbb Q G$-module homomorphism $\phi: \mathcal L \rightarrow \mathcal V$. Let $e_1=(1,0, \ldots, 0), \ldots, e_d=(0,\ldots,1)$ be the central idempotents of $(\mathbb Q G)^{(d)}$ and for $i=1, \ldots d$ write $e_i=x_i+y_i$ where $x_i \in \mathcal X, y_i \in \mathcal L$.  

Let $a \in \mathcal L$ and write $a=(l_1, \ldots, l_d)$ where $l_i \in \mathbb Q G$ so that $l_ie_i=ae_i$.
Thus \[ a=a (e_1+ \cdots +e_d)=\sum_{i=1}^d l_ie_i= \sum_{i=1}^d l_ix_i+ \sum_{i=1}^d l_iy_i.\]
Note that  \[ a \in \mathcal L , \quad \sum_{i=1}^d l_iy_i \in \mathcal L, \quad  \sum_{i=1}^d l_ix_i \in \mathcal X.\] Since $\mathcal X \cap \mathcal L=\{0\}$ this implies that $\sum_{i=1}^d l_ix_i=0$ and $a= \sum_{i=1}^d l_i y_i$. In particular $\phi(a)=\sum_{i=1}^d l_i \phi(y_i)$.

 We choose and fix $a \in \mathcal L$ such that $\phi(a) \not =0$.  For any integer $m >0$ we then have 
 \[  m^2 \phi(a)= \sum_{i=1}^d ml_i  \cdot m\phi(y_i) \not =0 .\]
 Now since $V$ and $\mathbb Z G$ are $\mathbb Z$-lattices in $\mathcal V$ and $\mathbb Q G$ respectively, we can pick an integer $m$ large enough so that $m\phi(y_i) \in V$ and $ml_i \in \mathbb Z G$ for all $i=1, \ldots, d$. We  set $h_i=m\phi(y_i)$ and $s_i=ml_i$ and record our conclusion as follows.

\begin{prop} There are elements $h_1, \ldots, h_d \in V$ and $(s_1, \ldots, s_d) \in L$ such that $z:=\sum_{i=1}^d s_ih_i \not = 0$.
\end{prop}

To finish the proof of Lemma \ref{module} let $c$ be any positive integer. Let $h_i,s_i$ and $z \in V$ be as above. Since $\cap_{i=1}^\infty  p^iV=\{0\}$ there is some integer $j$ such that $z \not \in p^jV$. Take $n> j+c$ and let $K_n$ be the submodule of $M$ generated by $\{ m_i+p^nh_i \ | \ i=1, \ldots, d\}$. We claim that the order of the element $z+K_n$ in the abelian group $M/K_n$ is a finite power of $p$ and at least $p^c$. On one hand as $L \subset \mathcal L = \ker f$ it follows that $\sum_{i=1}^ds_im_i=0$ and hence $\sum_{i=1}^d s_i(m_i+p^nh_i)=p^nz \in K_n$ and so the order of $z+K_n$ in $M/K_n$ is finite and divides $p^n$. Suppose the order of $z+K_n$ is less that $p^c$. This means that $p^kz \in K_n$ for some $k<c$. Recall the direct sum $U \oplus V \leq M$ and note that by definition $K_n \leq U \oplus p^nV$ since $h_i \in V$ for all $i=1, \ldots, d$. Now $p^kz \in K_n \cap V \leq p^nV$ and since $V$ is $\mathbb Z$-torsion free we obtain $z \in p^{n-k}V$. We have assumed that $n>j+c$ while $k<c$ and hence $z \in p^{n-k}V \subset p^jV$ in contradiction to the choice of $j$. Therefore the order of $z+K_n$ is at least $p^c$ and hence $t_p(M/K_n)\geq p^c$. The lemma follows.
$\square$

We say that a group $G$ is $p$-large if $G$ has a normal subgroup $N$ of index $p^k$ for some integer $k$ such that $N$ has a nonabelian free group as a homomorphic image.

The following result is an easy consequence of the main theorem of \cite{BT} and the work of Lackenby \cite{L} and Schlage-Puchta \cite{JCP}.
\begin{thm} \label{plarge} Let $F$ be the free group on generators $\{x,y\}$ and let $p$ be a prime. Let $n \in \mathbb N$ and let $a_j,b_j$, $j=1, \ldots, n$ be positive integers such that 
\[ \sum_{i=1}^n (p^{-a_j}+p^{-b_j}) <1.\]
Then for any elements $u_i,v_i\in F$ ($i=1,\ldots, n$), the group with presentation \[ G=\langle x,y \ | \ u_1^{p^{a_1}} \cdots u_n^{p^{a_n}}=1, \ v_1^{p^{b_1}} \cdots v_n^{p^{b_n}}=1 \rangle \]  is $p$-large.
\end{thm}
\textbf{Proof:} Let $H$ be the group with presentation 
\[ \langle x,y \ | \ u_1^{p^{a_1}}=v_1^{p^{b_1}}= \cdots = u_n^{p^{a_n}}=v_n^{p^{b_n}}=1  \rangle.\] Then $H$ is a homomorphic image of $G$. The presentation of $H$ above has $p$-deficiency (as defined in \cite{JCP}) greater than 0. By the main theorem of \cite{BT} the group $H$ is $p$-large and therefore so is $G$. $\square$

Before the last two results in this section we collect some basic properties of dimension of $p$-adic analytic groups. Their proofs can be found in \cite{DDMS}, Chapter 4.

\begin{prop} \label{dim} Let $G$ be a finitely generated $p$-adic analytic group. Then
\begin{itemize}
\item $\dim G= \dim H$ for every open subgroup of $G$.

\item $\dim G=0$ if and only if $G$ is a finite group.
 
\item Let $H$ be a closed normal subgroup of $G$. Then $\dim G= \dim H+ \dim G/H$. In particular $\dim G= \dim G/H$ if and only if $H$ is finite.
\end{itemize}
\end{prop}

\begin{prop} \label{G-1} Let $G$ be a finite group and let $A$ be an abelian pro-$p$ group which is a finitely generated $\mathbb Z_pG$-module. Then \[ t \left (\frac{ A}{(G-1)A} \right )\leq t(A) |G|^{\dim A}.\]
\end{prop}
\textbf{Proof:} Let $T$ be the torsion subgroup of $A$. Now \[ \frac{A}{ T +(G-1)A}= \frac{A/T}{(G-1) \cdot (A/T)}\] while $\frac{T+ (G-1)A}{(G-1) A}$ is an image of $T$ and hence has size at most $|T|=t(A)$. Therefore \[ t \left ( \frac{A}{(G-1)A} \right )=  \left | \frac{T+ (G-1)A}{(G-1) A} \right | t \left (\frac{A/T}{(G-1)(A/T)}\right) \leq t(A) t \left (\frac{A/T}{(G-1)(A/T)}\right) \] and hence it is enough to prove the proposition with $A/T$ replacing $A$. 

From now on we assume that $A$ is torsion free and hence a free $\mathbb Z_p$-module or rank $r:= \dim A$.
Let $W=A_G:= \{ a \in A \ | \ ga=a \ \forall g \in G\}$. 

We claim that $W \cap (G-1)A= \{0\}$. Define $z:= \sum_{g \in G} g \in \mathbb Z_p G$ and take $a \in (G-1) A \cap W$. Then $za=|G|a$ since $a \in W=A_G$ and on the other hand $za=0$ since $a \in (G-1)A$ thus $|G|a=0$ and so $a=0$ since $A$ is torsion free. The claim is proved. 

Suppose $b\in A$. Then $|G|b= zb- \sum_{ g \in G} (g-1)b$ and since $zb \in W$ and $b \in A$ was arbitrary we deduce that $|G|A \leq W \oplus (G-1)A$. In particular $A/(W \oplus (G-1)A) \leq |G|^r$, Therefore $A/(G-1)A$ has a torsion free subgroup isomorphic to $W$ of index at most $|G|^r$ and so $t(A/(G-1)A) \leq |G|^r$. The proposition follows.
$\square$

\begin{prop} \label{AB} Let $A$ be a finitely generated $p$-adic analytic group and let $B$ be an open normal pro-$p$ subgroup of $A$. 

1. If $\dim A^{ab}= \dim B^{ab}$ then $t(A^{ab})\leq t(B^{ab}) |A:B|$.

2. In all cases $t(A^{ab}) \leq t(B^{ab}) |A:B|^{\dim B^{ab}+1}.$
\end{prop}

\textbf{Proof:} Since $[B,B] \leq [A,A]$ we can factor out $[B,B]$ and from now on assume that $B$ is abelian. 

Part 1. The equality $\dim A= \dim B= \dim B^{ab}= \dim A^{ab}$ together with $\dim A= \dim A^{ab}+ \dim [A,A]$ this implies that $[A,A]$ is a finite group. Let $B= T \oplus E$ where $T$ is the torsion subgroup of $B$ and $E$ is torsion free. Note that $|A:E|= t(B) |A:B|$. Since $[A,A]$ is finite we have $[A,A] \cap E = \{1\}$ and in particular $E [A,A]$ is a torsion free subgroup of $A^{ab}$ which has index at most $|A:E| =t(B)|A:B|$. Therefore $t(A^{ab}) \leq t(B) |A:B|$ from Proposition \ref{el} and part 1 is proved. \medskip
 
Part 2. The abelian group $B$ is a finitely generated $\mathbb Z_p(A/B)$-module. Let $D$ be the torsion subgroup of $B/[A,B]$. By Proposition \ref{G-1} we have \[ |D| \leq t(B) |A:B|^{\dim B}.\] The group $B/[A,B]$ is a central subgroup of finite index in $A/[A,B]$ and therefore by Schur's theorem $[A,A]/[A,B]$ is finite. In particular the surjective homomorphism 
 
 \[ \pi: \frac{B}{[A,B]} \rightarrow \frac{B}{B \cap [A,A]}, \quad \pi(b\ [A,B])=b\ (B \cap [A,A]) \quad \forall b \in B  \] has kernel $(B \cap [A,A])/[A,B]$ which is a finite group. It follows that $t(\textrm{Im}\  \pi)= |D|/|\ker \pi|$ and therefore 
   \begin{equation}\label{D}  t(\textrm{Im}\ \pi)=t \left( \frac{B} {B \cap [A,A]}\right) \leq |D| \leq  t(B) |A:B|^{\dim B}. \end{equation} Observe that $B[A,A]/[A,A]$ is a subgroup of index at most $|A:B|$ in $A^{ab}$ and hence by Proposition \ref{el} \[t(A^{ab}) \leq |A:B| \ t\left (\frac{B[A,A]}{[A,A]} \right ).\] Together with the isomorphism $ \frac{B[A,A]}{[A,A]} \simeq \frac{B}{B \cap [A,A]}$ and (\ref{D}) this proves part 2.
$\square$
\section{Proof of Theorem \ref{main}}

Let $F$ be the free group on two generators $x,y$ and let $f: \mathbb N \rightarrow \mathbb N$ be any function.
\begin{prop}\label{seq}
There exist
\begin{itemize}
\item A sequence $F=F_0>F_1>F_2 \cdots$ of normal subgroups of $F$, with $|F:F_i|=p^{q_i}$,  \item Two sequences of $(r_i)_{i=0}^\infty,(w_i)_{i=0}^\infty$ of elements in $\cap_{j=1}^\infty F_j$, 
\item Two associated sequences of groups $(N_i)_{i=0}^\infty, (H_i)_{i=0}^\infty$ with $N_i$ defined as the normal subgroup generated by $\{r_i,w_i\}$ in $F$ and $H_i:=F_i/N_i$, \end{itemize}
such that the following properties hold:

P1: $t_p(H_i^{ab})>p^{f(q_i)}$ for each $i \geq 1$.

P2. For each $i \geq 1$ there is an integer $a_i \geq 3$ and elements $u_i, v_i \in F_{i}$ such that $r_{i}=r_{i-1} u_i^{p^{a_i}}$ and $w_{i}=w_{i-1} v_i^{p^{a_i}}$. Moreover for $i>1$ we have 
\[  a_i > \max \{ a_{i-1}, \log_p t_p (H_{i-1}^{ab})  \}.\]
   \end{prop}
   
 \textbf{Proof:} We set $F_0=F$ and $r_0=w_0=1$.
   We will construct $F_i,r_i,w_i$ by induction on $i$.
   
Suppose $i \geq 1$ and $F_j,r_j,w_j$ have been chosen for $j=1, \ldots, i-1$ subject to P1 and P2.
We show the existence of $F_i,r_i$ and $w_i$ as follows:

 From P2 we deduce $r_{i-1}=u_1^{p^{a_1}} \cdots u_{i-1}^{p^{a_{i-1}}}$ and 
 $w_{i-1}=v_1^{p^{a_1}} \cdots v_{i-1}^{p^{a_{i-1}}}$.  In addition $3 \leq a_1 <a_2 < \cdots $
gives  $a_j \geq 2+j$ for $j=1, \ldots i-1$. Therefore
\begin{equation}\label{deficiency}  \sum_{j=1}^{i-1} (p^{-a_j}+p^{-a_j})< 2 \sum_{j=1}^\infty p^{-2-j}=\frac{2}{p^3(1-p^{-1})}=\frac{2}{p^2(p-1)} < 1.\end{equation}  
We claim that $G_{i-1}=F/N_{i-1}$ contains a normal subgroup $S$ of $p$-power index with infinite abelianization and such that $S<H_{i-1}=F_{i-1}/N_{i-1}$.
 Theorem \ref{plarge} together with (\ref{deficiency}) gives that the group $G_{i-1} \simeq  \langle x,y \ | \ r_{i-1}=w_{i-1}=1\rangle$ is $p$-large. In particular $G_{i-1}$ has a normal subgroup of $p$-power index $L$ such that $L^{ab}$ is infinite. Replacing $L$ with a smaller $G_{i-1}$-normal subgroup of $p$-power index and contaning $[L,L]$, we may assume in addition that $|G_{i-1}:L|>|G_{i-1}:H_{i-1}|$. Therefore $L \not \geq H_{i-1}$. Let $S=L \cap H_{i-1}$, this is a normal subgroup of $p$-power index in $G_{i-1}$ and $S< H_{i-1}$. Since $S[L,L]/[L,L]$ has finite index in $L/[L,L]=L^{ab}$ we deduce that $S^{ab}$ is infinite, proving the claim.
 
 Let $F_i$ be the subgroup of $F$ such that $S=F_i/N_{i-1}$ and let $|F:F_i|=p^{q_i}$.
 Note that $S^{ab}=F_i/[F_i,F_i]N_{i-1}$ and hence $|F_i: [F_i,F_i]N_{i-1}|$ is infinite.
 
 It remains to specify $u_i^{a_i}$ and $v_i^{a_i}$. Let $G:= F/F_i$ and let $M:=F_i^{ab}$ be the relation module of the associated presentation \[ 1 \rightarrow F_i \rightarrow F \rightarrow G \rightarrow 1\] of $G$. Further let $m_1=r_{i-1}[F_i,F_i]$ and $m_2= w_{i-1}[F_i,F_i]$. The $\mathbb Z G$-submodule generated by $m_1,m_2$ in $M$ equals $\frac{N_{i-1}[F_i,F_i]}{[F_i,F_i]}$ and has infinite index in $M$ by the choice of $S$.
 The relation module $M$ embeds in $(\mathbb Z G)^{(2)}$ (see \cite{B} \S II, Proposition 5.4) and in particular $M$ is $\mathbb Z$-torsion free. Moreover $\mathbb Q \otimes_\mathbb Z M$ is a submodule of $(\mathbb Q G)^{(2)}$ and hence also a direct factor of it. Thus $\mathbb Q \otimes_\mathbb Z M$ is $2$-generated as $\mathbb Q G$-module. All conditions of Lemma \ref{module} are satisfied. Let $h_1,h_2 \in M$ and $K_n$ be as in the conclusion of the Lemma and let $N \in \mathbb N$ be such that for all $n>N$ we have $t_p(M/K_n)>p^{f(q_i)}$. Choose elements $u_i, v_i \in F_i$ such that $h_1=u_i[F_i,F_i]$ and $h_2=v_i [F_i,F_i]$. Choose an integer $a_i$ such that \[ a_i>\max\{N,a_{i-1}, \log_p t_p (H_{i-1}^{ab})\}.\] Let $r_i=r_{i-1}u_i^{p^{a_i}}$ and $w_i=w_{i-1} v_i^{p^{a_i}}$ so that $r_i[F_i,F_i]=m_1+p^{a_i}h_1$ and $w_i[F_i,F_i]=m_2+p^{a_i}h_2$.
 It follows that $K_{a_i}=N_i[F_i,F_i]$ and therefore 
 \[ t_p(H_i^{ab})= t_p\left (\frac{F_i}{[F_i,F_i]N_i} \right )=t_p \left (\frac{M}{K_{a_i}}\right )>p^{f(q_i)}. \]
 
 Hence $F_i,r_i,w_i$ satisfy P1 and P2 as required. The induction step is complete. $\square$
 
\medskip

For the rest of this section let $F_i,r_i=u_1^{p^{a_1}} \cdots u_{i}^{p^{a_{i}}}$ and $w_i=v_1^{p^{a_1}} \cdots v_{i}^{p^{a_{i}}}$ be as given in Proposition \ref{seq}. In particular $N_i$ is the normal closure of $\{r_i, w_i\}$ in $F$, and $H_i=F_i/N_i$.

Let $\hat F$ be the pro-$p$ completion of $F$ and denote by $\hat F_i$ the closure of $F_i$ in $\hat F$. It is well known that $\hat F_i$ is the pro-$p$ completion of $F_i$.

From $P2$ it follows that $r_j \in r_i F^{p^{a_{i+1}}}$ for all $1 \leq i<j$ and in particular $r_1,r_2, \ldots$ is a Cauchy sequence in $\hat F$. Hence there is a unique element $r_\infty \in \hat F$ such that $\lim_{i \rightarrow \infty} r_i = r_\infty$ and
similarly we define $w_\infty \in \hat F$ as $w_\infty=\lim_{i \rightarrow \infty} w_i$. 

Observe that P2 imples  

\begin{equation} \label{[n]} r_i \equiv  r_\infty  \ \textrm{ and } \  w_i \equiv  w_\infty \quad \textrm{mod} \ \hat F_i^{p^{a_{i+1}}} \quad \forall i \in \mathbb N. \end{equation}

Let $Z_\infty$ be the normal closure of $\{r_\infty,w_\infty\}$ in $\hat F$, that is the smallest normal closed subgroup of $\hat F$ containing $\{r_\infty,w_\infty\}$. Further, for $i=1, 2, \ldots $ let $Z_i$ be the normal closure of $r_i,w_i$ in $\hat F$. Define $\Gamma= \hat F/Z_\infty$ and $\Gamma_i=\hat F_i/Z_\infty$. Then $|\Gamma:\Gamma_i|=|\hat F:\hat F_i|=|F:F_i|=p^{q_i}$

Recall that for any group we denote $G[n]= G/G^n$. The congruences (\ref{[n]}) give that \[ Z_i \hat F_i^{p^{a_{i+1}}}  [\hat F_i, \hat F_i]=Z_{\infty} \hat F_i^{p^{a_{i+1}}} [\hat F_i, \hat F_i] \] and  therefore \begin{equation} \label{Gam} \Gamma_i^{ab}[p^{a_{i+1}}]=\hat H_i^{ab}[p^{a_{i+1}}], \end{equation} where $\hat H_i=\hat F_i/Z_i$. 

Choose and fix a positive integer $i$ for the rest of this section. Let $\bar N_i$ be the closure of $N_i$ in $\hat F_i$. Since $\hat F_i$ is closed in $\hat F$ we note that $\bar N_i$ equals the closure of $N_i$ in $\hat F$.

We claim that $Z_i=\bar N_i$. Indeed $\bar N_i$ is a closed subgroup of $\hat F$ which contains $r_i,w_i$ and which is normalized by $F$ and hence by all of $\hat F$. Therefore $Z_i \leq \bar N_i$ by the definition of $Z_i$. 
On the other hand $Z_i$ is normal in $\hat F$ and so $Z_i \cap F$ is normal in $F$ and contains $r_i,w_i$ and hence their normal closure $N_i$ in $F$. Therefore $Z_i\geq N_i$ and so $Z _i\geq \bar N_i$ since $Z_i$ is closed by definition.  Hence $Z_i=\bar N_i$ proving the claim. 

Since $\hat F_i$ is the pro-$p$ completion of $F_i$ it follows that $\hat H_i=\hat F_i/Z_i=\hat F_i/ \bar N_i$ is the pro-$p$ completion of $H_i=F_i/N_i$.

In particular $\hat H_i^{ab}$ is the pro-$p$ completion of $H_i^{ab}$ and therefore $t(\hat H_i^{ab})= t_p(H_i^{ab})$. We have $a_{i+1}> \log_p t(\hat H_i^{ab})$ and by Lemma \ref{L1} applied to the abelian pro-$p$ groups $\hat H_i^{ab}$ and $\Gamma_i^{ab}$ together with (\ref{Gam}) above we deduce \[ t(\Gamma_i^{ab}) \geq t(\hat H_i^{ab})\quad \forall i \geq 1.\] On the other hand $t(\hat H _i^{ab})=t_p(H_i^{ab})>p^{f(q_i)}$ by property P1. Combining the last two inequalities we deduce $t(\Gamma_i^{ab}) > p^{f(q_i)}$ and Theorem \ref{main} is proved. \medskip

We note that $\Gamma$ maps onto the group $\Delta :=\hat F/W$, where $W$ is the normal closure of $\{ u_i^{p^{a_i}}, v_i^{p^{a_i}} \ | \ i \in \mathbb N\}$. The presentation $\hat F/W$ defining $\Delta$ has positive $p$-deficiency in view of (\ref{deficiency}). It is easy to show that $\Delta$ is also a Golod-Shafarevich group.

\section{Proof of Theorem \ref{padic}}
\subsection{When $G$ is uniform.}
Suppose first that $G$ is a uniformly powerful pro-$p$ group. Let $\mathcal G=(G,+, (,))$ be the Lie ring structure on $G$ as defined in \cite{DDMS}, Chapter 4.5. We will denote the Lie bracket of two elements $x, y \in \mathcal G$ by $(x,y)$ and let  $x+y$ be their sum in the additive group of $\mathcal G$. Since $G$ is uniformly powerful, both $G$ and $(\mathcal G, +)$ are torsion free groups. Moreover $[G,G]=(\mathcal G, \mathcal G)$ by Theorem B of \cite{sanches}. We have $x+ y \equiv xy $ mod $[G,G]$ from the Baker-Campbell-Hausdorff formula, see Theorem 4.5 of \cite{sanches} and therefore the identity map induces an isomorphism between $G^{ab}=G/[G,G]$ and the additive group of $\mathcal G/(\mathcal G, \mathcal G)$. In particular \[t(G^{ab})= t \left (\frac{\mathcal G}{(\mathcal G, \mathcal G)}\right).\]

For $n \in \mathbb N$ we write $G_n= G^{p^n}= \{ g^{p^n}  \ | \ g \in G \}$. Then $G_n$ is also a uniformly powerful pro-$p$ group with Lie structure $\mathcal G_n=p^n \mathcal G$.
 Since the additive group of $\mathcal G$ is torsion free the map $x \mapsto p^n x$ induces an isomorphism $\mathcal G/p^n (\mathcal G, \mathcal G) \rightarrow p^n\mathcal G/(p^n \mathcal G, p^n \mathcal G) $. In turn $(\mathcal G,\mathcal G)/p^n (\mathcal G, \mathcal G)$ is a finite group of size at most $p^{n \dim \mathcal G}= p^{n \dim G}$ and therefore
\[ t(G_n^{ab})= t \left ( \frac{p^n\mathcal G}{(p^n\mathcal G, p^n\mathcal G)} \right) = t \left ( \frac{\mathcal G}{p^{n}(\mathcal G, \mathcal G)} \right)= t\left (\frac{\mathcal G}{(\mathcal G, \mathcal G)}\right ) \left |\frac{(\mathcal G,\mathcal G)}{p^{n}(\mathcal G, \mathcal G)} \right| \leq a p^{n \dim G}, \] where $a=t(G^{ab})$.

We also note that $[G_n,G_n]= (\mathcal G_n, \mathcal G_n)= p^{2n}(\mathcal G, \mathcal G)$ and $[G,G]= (\mathcal G, \mathcal G)$ imply \[ \dim [G,G]= \dim [G_n,G_n].\] 

Suppose $H$ is a normal open subgroup of $G$ and let $|G:H|=p^n$. Then $G_n=G^{p^n} \leq H$. Observe that $[G_n,G_n] \leq [H,H] \leq [G,G]$ and therefore $\dim [H,H]=\dim[G_n,G_n]= \dim [G,G]$. Together with $\dim H=\dim G=\dim G_n$ we obtain  $\dim H^{ab}= \dim G^{ab}= \dim G_n^{ab}$. Proposition  \ref{AB} Part 1 with $A=H$ and $B=G_n$ gives 
\begin{equation}\label{uniform}  t(H^{ab}) \leq t(G^{ab}_n) |H: G_n| \leq a p^{2n \dim G}=a |G:H|^{2 \dim G}.\end{equation}
since $|H:G_n| \leq |G:G_n|=p^{n \dim G}$ and $t(G_n^{ab}) \leq a p^{n\dim G}$.

\subsection{The general case.}
Finally consider the case when $G$ is an arbitrary $p$-adic analytic pro-$p$ group. Choose and fix $G_0$, a uniform open normal pro-$p$ subgroup of $G$. Let $H$ be any open normal subgroup of $G$ and define  $H_0=G_0 \cap H$. Then $|H:H_0|=|HG_0:G_0| \leq |G:G_0|$ while $|G_0:H_0|= |G_0H:H| \leq |G:H|$. The inequality (\ref{uniform}) applied to the open subgroup $H_0$ of $G_0$ gives \[ t(H_0^{ab}) \leq a_0 |G_0:H_0|^{2 \dim G} \leq a_0|G:H|^{2 \dim G},\] where $a_0=t(G_0^{ab})$. 

Now apply Proposition \ref{AB} Part 2 with $A=H$ and $B=H_0$ noting that $\dim H_0^{ab} \leq \dim G$. This gives $t(H^{ab}) \leq t(H_0^{ab}) |H:H_0|^{\dim G+1}$ and together with the bound for $t(H_0^{ab})$ we obtain
\[ t(H^{ab}) \leq a_0|H:H_0|^{\dim G+1} |G:H|^{2 \dim G} \leq  b |G:H|^{2 \dim G},  \] where $b=a_0 |G:G_0|^{\dim G+1}$ is a constant depening on $G$ and $G_0$ but not on $H$. Theorem \ref{padic} is proved. $\square$


\begin{thebibliography}{99}
\bibitem{AGN} 
M. Abert, T. Gelander, N. Nikolov, Rank, combinatorial cost and homology torsion growth in higher rank lattices, Duke Math. J., vol. 166 (15) (2016) pp. 2925-2964. 

\bibitem{BSV} N. Bergeron, M. H. Sengun, A. Venkatesh. Torsion homology growth and cycle
complexity of arithmetic manifolds. Duke Math. J., 165,  (2016), pp. 1629-1693.

\bibitem{B} K. S. Brown, Cohomology of Groups, Graduate texts in Mathematics 87, Springer, 1982.

\bibitem{BT}  J. O. Button, A. Thillaisundaram, Applications of $p$-deficiency and $p$-largeness, International J. Algebra and Computation Vol. 21, No. 04,  (2011) pp. 547-574.

\bibitem{DDMS} J. Dixon, M. du Sautoy, A. Mann, D. Segal, Analytic pro-$p$ groups, 2nd ed. Cambridge University Press, 1999.
 
\bibitem{sanches} J. Gonzales-Sanches, On $p$-saturable groups, J. Algebra 315 (2007) pp. 809-823.

\bibitem{KKN} A. Kar, P. Kropholler, N. Nikolov, On growth of homology torsion in amenable groups, Math. Proc. Cambridge Phil. Soc., 162 (2), pp. 337-351.

\bibitem{L} M. Lackenby, Detecting large groups, J. Algebra, Vol. 324 (2010), pp. 2636-2657.

\bibitem{JCP} J. C. Schlage-Puchta, A $p$-group with positive rank gradient, J. Group Theory 15 (2012), pp. 261-270.

\end{thebibliography}
\end{document}